\documentclass[12pt]{article}
\usepackage{amsmath,amsthm,amssymb,epsfig,graphics,subfigure,color}
\usepackage[english]{babel}
\input amssym.def
\input amssym.tex
\date{}

\newtheorem{theo}{Theorem}
\newtheorem{lem}{Lemma}
\newtheorem{cor}{Corollary}

\newcommand\Lip{\mathop{\mathrm{Lip}}\nolimits}

\begin{document}

\title{{\bf On a problem of Eidelheit from {\tt The Scottish Book}
concerning absolutely continuous functions}}

\author {\large Lech Maligranda\thanks{Research partially
supported by the Swedish Research Council (VR) grant
621-2008-5058},\, Volodymyr V. Mykhaylyuk\\
{\small and} \large  Anatolij Plichko}

\date{}

\maketitle

\renewcommand{\thefootnote}{\fnsymbol{footnote}}


\vspace{-5mm}

\hspace{50mm}{\it Dedicated to the memory of M. Eidelheit (1910-1943)}

\hspace{60mm}{\it on the occasion of $100^{\small th}$ years of his birth}

\begin{abstract}.

\vspace{-5mm}

\noindent {\footnotesize A negative solution of Problem 188
posed by Max Eidelheit in {\it the Scottish Book} concerning
superpositions of separately absolutely continuous functions is
presented. We discuss here his and some related problems which have
also negative solutions. Finally, we give an explanation of such negative
answers from the "embeddings of Banach spaces" point of view.}
\end{abstract}

\section{Introduction}

\noindent
There are several equivalent definitions of the concept absolute
continuity. The notion and the term of {\it absolutely continuous}
was introduced in 1905 by G. Vitali \cite{Vitali}. Let $I=[a,
b]\subset \mathbb R$ and $f:I\to\mathbb R$. The function $f$ is called
{\it absolutely continuous on $I$} if for every $\varepsilon>0$
there exists a $\delta>0$ such that for any $a\le a_1<b_1\le a_2<
b_2\le \dots \le a_n< b_n\le b$ the condition
$\sum_{k=1}^{n}(b_k-a_k)<\delta$ implies that
$\sum_{k=1}^{n}|f(b_k)-f(a_k)|<\varepsilon$ (cf. Natanson
\cite{Natanson}, p. 243). Also we can say that for every
$\varepsilon>0$ there exists a $\delta>0$ such that for any finite
collection of mutually disjoint intervals $I_k=(a_k,b_k)\subset
I\;(k =1,2,\ldots,n)$ we have that $\sum_{k=1}^{n}|I_k|<\delta$ implies
$\sum_{k=1}^{n}|f(I_k)|<\varepsilon$. The requirement that the
open intervals $I_k$ are disjoint is sometimes stated by saying
that the corresponding closed intervals $\{[a_k,b_k]\}_{k=1}^n$
must be {\it nonoverlapping}, that is, their interiors are
disjoint. Note that since the number $n\in \mathbb N$ is arbitrary, we
can also take $n=\infty$, that is, replace finite sums by series.

It is obvious that an absolutely continuous function is continuous
and it is easy to show that it is also of bounded variation. The
classical Banach-Zarecki{\u \i} theorem states that a function
$f:I\to\mathbb R$ is absolutely continuous if and only if it is
continuous, is of bounded variation and has the Luzin (N)
property, that is, maps null sets into null sets
(cf. \cite[Th. 7.11]{BBT}, \cite[p. 250]{Natanson} and
\cite[p. 146]{Rudin}). Of course, every {\it Lipschitz function}
on $I$, that is, any function $f:I \to \mathbb R$ satisfying the condition
that there exists a constant $C>0$ such that $|f(x)-f(y)| \leq C
|x-y|$ for all $x, y \in I$, is absolutely continuous on $I$. All
Lipschitz functions on $I$ we denote by $\Lip_1(I)$. One of
the equivalent norms in $\Lip_1(I)$ is defined by
$$
\|f\|=|m_f|+\sup_{x,y\in I,\;x\ne y}\frac{|f(x)-f(y)|}{|x-y|}, ~~ {\rm where} ~~ m_f=\int_If(x)\,dx.
$$

The fundamental theorem of calculus for absolutely continuous
functions (cf. \cite[pp. 253--255]{Natanson}, \cite[pp.
197--198]{Rana} and \cite[pp. 148--149]{Rudin}) gives that a
function $f:I\to\mathbb R$ is absolutely continuous if and only if
$f$ is differentiable almost everywhere on $I$, the derivative $f'
\in L_1(I)$, i.e. is Lebesgue integrable, and $f(x)=f(a)+\int_a^x
f'(t)\,dt$ for every $x\in I$. On the other hand, if we put
instead of $f'\in L_1(I)$ a stronger assumption $f' \in
L_{\infty}(I)$ we obtain a characterization of Lipschitz functions
on $I$. Therefore, for an absolutely continuous function $f$, the
condition $\int_a^b|f'(x)|^p\,dx <\infty$ for each $p>1$ is a
natural weakening of the Lipschitz condition.\medskip

In what follows we consider only the segment $I=[a, b]=[0,1]$ and
the square $Q=I^2$.

\medskip
Max Eidelheit on October 27, 1940 wrote in {\it The Scottish Book}
the following problem concerning superposition of absolutely
continuous functions (cf. \cite{Mauldin}, Problem 188.1, p.
261$^{\small 1}$):
\footnotetext{$^{\small 1}$ In the original handwritten {\it
Scottish Book} in Polish language this problem has number 188 and
there was double numeration of the Problem 185 (one written by
Saks and the other one by Banach). In the English translation done
by Ulam in 1957 appeared instead double numeration of the Problem
188 (one by Sobolev, which originally has number 187, and the other
one by Eidelheit). This was probably the reason why in the Mauldin
edition of {\it The Scottish Book} \cite{Mauldin} we have the
numbers 188 on the Sobolev problem and 188.1 on the Eidelheit
problem. One can suppose that the integral conditions of the
Eidelheit problem are connected with Sobolev's visit to the {\it
Scottish Caf\'e} after which Eidelheit became more familiar with
Sobolev spaces.}
\medskip

\noindent{\bf Problem} (Eidelheit). {\it Let a function
$f:Q\to\mathbb R$ be absolutely continuous on every straight line
parallel to the axes of the coordinate system and let
$g_1,g_2:I\to I$ be absolutely continuous functions. Is the
function $f(g_1(t),g_2(t))$ also absolutely continuous? If not,
then perhaps this holds under the additional assumptions that
$\mathop{\int\!\!\int}\limits_{Q} |f'_x|^p\, dxdy<\infty$ and
$\mathop{\int\!\!\int}\limits_{Q}|f'_y|^p\, dxdy<\infty$, where
$p>1$?}

\medskip
There is no any comment to this problem in the book
\cite{Mauldin} on page 261.

\medskip
Note that there are several different meanings of the conditions
in the problem: $f$ can be absolutely continuous on every straight
line parallel to the axes or on almost every straight line parallel
to the axes, the integrals can be bounded for some $p$ or for every $p$
and the derivatives can exists everywhere or almost everywhere.

\medskip
Our intention here is to give some short historical comments to the
Eidelheit problem (as we will show in Theorem \ref{solution}, the
answer has been known to a great extent even before the problem was posed),
and present some variations and generalizations of known results
connected with this problem.

\medskip
It easy to see that the first part of Eidelheit's question has a
negative answer. Consider the Schwarz function
\begin{equation*}
f(x, y) = \left\{
\begin{array}{ll}
\frac{2xy}{x^2+y^2}, & ~~~\mbox{{\rm if} ~$x^2 + y^2 > 0$,} \\
0, & ~~~\mbox{{\rm if} ~$x = y = 0$.}\\
\end{array}
\right.  \label{def_phicomp}
\end{equation*}
The function $f$ is absolutely continuous in each variable since for
any fixed $y\in (0,1]$ we have that
$|f(x,y)-f(u,y)|\le\frac{2}{y}|x-u|$ for all $x,u\in I$ and
$f(x,0)-f(u,0)=0$. Similarly with fixed $x$. If we take the
functions $g_1(t)=g_2(t)=t$, then the superposition $f(g_1(t),
g_2(t))=f(t,t)$ is $2$ for $t\ne 0$ and $0$ for $t=0$ and, hence, it is
discontinuous at $t=0$, and so not absolutely continuous on $I$.
Note that the integrals of Eidelheit problem are unbounded for the
Schwartz function if $p>1$.
\smallskip

Also the second part of Eidelheit's problem has a negative answer,
which we will present in the next section. We even give a negative
answer to the diagonal case, that is, when $g_1(t) = g_2(t) = t$.

The paper is organized as follows: In Section 2 we show how to
obtain the answer to Eidelheit's question using a well known
Fichtenholz theorem. Then we receive two variable Fichtenholz
theorems concerning superposition of absolutely continuous functions as a
corollary of a general theorem on superpositions in Banach spaces.
Section~3 contains a counterexample to the diagonal version
of Eidelheit's problem. Finally, in Section 4 we give the
``embeddings of Banach spaces'' approach to this problem.

\section{Superposition of absolutely continuous functions}

We start with the question about the superposition of one variable
functions. As is well known, the functions $f(x)=x^{1/2}$ and
$g(x)= x^2\sin^2(\frac{1}{x})$, if $x>0$ and $=0$, if $x=0$ are
absolutely continuous on $I$ but their superposition $f\circ g$ is
not since it has infinite variation. We even have that
$g\in\Lip_1(I)$ since $|g'(x)|\le 4$ for all $x\in I$.

For the first time the existence of superposition of absolutely continuous
functions which is not absolutely continuous was noticed by
W. Wilkosz \cite[pp. 479-480]{Wilkosz} (cf. also \cite[Ex. 1]{Rakoczy})
and D. Jackson. More exactly, it was noted in \cite[p. 462]{VallePoussin}
that the proposition from \cite[p. 280]{Valle} on absolute continuity
of the superposition of absolutely continuous functions is not true
and that Jackson informed de la Vall\'ee Poussin about
this fact. In the Carath\'eodory book \cite[pp. 554--555]{Cara}
there is a construction of $g: I \rightarrow I$ such that its
variation is equal to $1$, but the variation of $\sqrt{g}$ is
infinite. The absolute continuity of the superposition of
absolutely continuous functions was investigated in detail
by G. Fichtenholz. Already in 1922 (cf. \cite{Fichtenholz22},
pp.~436--439 and \cite{Fichtenholz}, pp.~289--291) Fichtenholz
proved the following result, where we can see big difference
between nonoverlapping and overlapping intervals:

\medskip
\noindent{\bf Theorem A (Fichtenholz).} {\it Let $f:I\to\mathbb R$
be a function. Then the following conditions are equivalent:
\begin{itemize}
\item[$(i)$] $f$ is Lipschitz on $I$. \item[$(ii)$] For every
$\varepsilon>0$ there exists $\delta>0$ such that for any $0\le
a_k<b_k\le 1, (k =1,2,\ldots n)$ the condition
$\sum\nolimits_{k=1}^{n}(b_k-a_k)<\delta$ implies
$\sum\nolimits_{k=1}^{n}|f(b_k)-f(a_k)|<\varepsilon$.
\item[$(iii)$] For every absolutely continuous function $g:I\to I$
the superposition $f\circ g$ is absolutely continuous on $I$.
\item[$(iv)$] For every Lipschitz function $g:I\to I$ the
superposition $f\circ g$ is absolutely continuous on $I$.
\end{itemize}}

The Banach-Zarecki{\u \i} theorem indicates that a superposition of
two absolutely continuous functions can fail to be absolutely
continuous if and only if it is not of bounded variation since
both continuity and Luzin's condition (N) are preserved under
superposition. Therefore the question about superposition of functions
of bounded variation has the same answer as that about superposition
of absolutely continuous functions. From the Fichtenholz
characterization in Theorem A we can get a similar characterization
for functions of bounded variation (BV), which was done in 1981 by
Josephy \cite[Th. 2]{Josephy}: for $f\colon I\rightarrow I$ the
superposition $f\circ g\in\text{BV}$ for all $g\in\text{BV}$ if
and only if $f$ is a Lipschitz function on $I$.
\medskip

It is well-known since a long time that the absolutely continuous
function $f(x)=\int_0^x\ln t\,dt$ is not Lipschitz, because its
derivative $f'(x)=\ln x$ is unbounded, and that $\int_0^1|\ln
x|^p\,dx<\infty$ for every $p>1$. This fact together with Theorem A
answers the one variable version of Eidelheit's question.
A negative answer to one variable Eidelheit's question gives
automatically the negative answer to the corresponding two
(and $n$) variable question. More exactly, we have:

\begin{theo}\label{solution}
There exists a function $\varphi:Q \to\mathbb R$, absolutely
continuous in each variable, such that
$\mathop{\int\!\!\!\int}\limits_{Q}|\varphi'_x|^p \,dxdy<\infty$,
$\mathop{\int\!\!\!\int}\limits_{Q}|\varphi'_y|^p \,dxdy<\infty$
for every $p>1$, and Lipschitz functions $g_1,g_2$ in $I$ such
that the superposition $\varphi(g_1(t),g_2(t))$ is not absolutely
continuous.
\end{theo}

{\it Proof.} Indeed, let $f(x)=\int_0^x\ln t\,dt$ and $g$ be the
corresponding Lipschitz function from Theorem A. Put
$\varphi(x,y)=f(x)$ (the function $\varphi$ depends on $y$ only
formally) and $g_1=g_2=g$. Then $\varphi$ is absolutely continuous
in each variable and for each $p$ we have that
$$
\mathop{\int\!\!\!\int}_Q|\varphi'_x|^p\,dxdy=\int\limits_I|f'|^p\,dx<\infty\;,
\;\;\mathop{\int\!\!\!\int}_Q|\varphi'_y|^p\,dxdy=0,
$$
and the superposition $\varphi(g_1(t),g_2(t))=f(g(t))$ is not absolutely
continuous.  \qed
\bigskip

Naturally there arises a question on the validity of the two variables
(and $n$ variables) Fichtenholz theorem. We even prefer to obtain this
generalization as a corollary of a more general theorem on
superpositions in normed spaces.
\medskip

Let us recall that a mapping $f: X \rightarrow Y$
defined on a metric space $(X,d)$ with values in a metric space
$(Y,\rho)$ is called {\it Lipschitz on $X$} if there exists a
constant $C>0$ such that
$$
\rho(f(x),f(y))\le C\,d(x,y) ~~ \mbox{for all}~~ x,y \in X.
$$
A mapping $g:I\to X$ is called {\it absolutely continuous on
$I$} if for every $\varepsilon>0$ there is $\delta>0$ such that
for any $0\le a_1<b_1\le a_2< b_2\le \dots\le a_n< b_n\le 1$ the
condition $\sum_{k=1}^{n}(b_k-a_k)<\delta$ implies that
$\sum_{k=1}^{n}d(g(b_k),g(a_k))<\varepsilon$.

\begin{theo}
Let $P$ be a convex set in a normed space $X$ and $f:P\to\mathbb R$ be
a function. Then the following conditions are equivalent:
\begin{itemize}
\item[$(i)$] For every compact set $K\subset P$ the restriction
$f|_K$ is Lipschitz on $K$.

\item[$(ii)$] For every absolutely continuous mapping $g:I\to P$
the superposition $h = f\circ g$ is absolutely continuous on $I$.

\item[$(iii)$] For every Lipschitz mapping $g:I\to P$ the
superposition $h=f\circ g$ is absolutely continuous on $I$.
\end{itemize}
\end{theo}

{\it Proof.} $(i)\Rightarrow (ii)$. Let $g:I\to P$ be absolutely
continuous. Since the set $K=g(I)$ is compact as a continuous
image of the compact set $I$, the function $f\circ g$ is absolutely
continuous on $I$ as a superposition of Lipschitz and absolutely
continuous mappings.\medskip

The implication $(ii)\Rightarrow (iii)$ is obvious.\medskip

$(iii)\Rightarrow (i)$. Suppose that there exists a compact set
$K\subset P$ for which the restriction $f|_K$ is not Lipschitz.
Then there are points $x_n,y_n\in P$, $n=1,2,\dots,$ such that for
every $n \in \mathbb N$ we have
$$
\frac{|f(x_n)-f(y_n)|}{d_n}\ge 2n^3\max_{x\in K}|f(x)|\,,
$$
where $d_n=\|x_n-y_n\|$. Hence, $d_n\leq\frac{1}{n^3}$
for every $n$. Since $K$ is compact, without loss of generality,
one can assume that
\medskip

\hspace{1,7 cm} $\sum_{n=1}^{\infty}\|x_n-x_{n+1}\|<\infty$, in particular,
there exists $x_0=\lim_{n\to\infty}x_n$.
\medskip

For every $n\in \mathbb N$ put $k_n=[\frac{1}{n^2d_n}]$. Note that
$\frac{1}{n^2}-d_n<k_nd_n\le \frac{1}{n^2}$. Hence, $k_nd_n\sim
\frac{1}{n^2}$ and $\sum_{n=1}^{\infty}
k_nd_n<\infty$.
\smallskip

Define a sequence of segments $[a_n,b_n]$ recursively as follows: put
$a_1=0$ and for each $n>0$,
$$
b_n=a_n+2k_nd_n\;,\;\;\; a_{n+1}=b_n+\|x_n-x_{n+1}\|\,.
$$
Let $b=\lim_nb_n$. Then $b<\infty$, since $\sum_{n=1}^{\infty} k_nd_n<\infty$ and $\sum_{n=1}^{\infty}\|x_n-x_{n+1}\|<\infty$.
By definition
$$
0=a_1<b_1\le a_2<b_2\le\dots
$$
Moreover,
$$
[0,b]=\bigcup_{n=1}^{\infty}([a_n,b_n]\cup[b_n,a_{n+1}]) \cup\{b\}.
$$

Let us construct a Lipschitz mapping $g:[0,b]\to P$ such that
the composition function $h:[0,b]\to\mathbb R$ given by $\;h(x)=f(g(x))$ is not
absolutely continuous.$^{\small 2}$
\footnotetext{$^{\small 2}$ In this construction we use the interval $I=[0,b]$ instead of
$I=[0,1]$, but this is not essential.}

\medskip
\noindent
We define the mapping $g$, on each segment $[a_n,b_n]$,
as follows:
\medskip

(1) $g(a_n+2id_n)=x_n\;$ for $\;0\le i\le k_n$ and $\;g(a_n+(2i-1)d_n)=y_n\;$
for $\;1\le i\le k_n$;
\smallskip

(2) $g\;$ is linear on every segment $\;[a_n+(j-1)d_n,
a_n+jd_n]\;$ for $\;1\le j\le 2k_n$.
\medskip

\noindent
The length of each segment $[a_n+(j-1)d_n,
a_n+jd_n]$ is equal to $d_n=\|x_n-y_n\|$, so $g$ is Lipschitz
with the constant $C=1$ on each such segment (hence on the whole
interval $[a_n,b_n]$).

Define $g$ to be linear on each segment
$[b_n,a_{n+1}]$. Since, by $(1)$, $g(b_n)=x_n\,$,
$\;g(a_{n+1})=x_{n+1}$, and
$\|x_n-x_{n+1}\|=a_{n+1}-b_n$, the mapping $g$ is
Lipschitz with the constant $C=1$ on the segment
$[b_n,a_{n+1}]$. Finally, $g(b)=x_0$.

Since $g$ is Lipschitz with the constant $C=1$ on $[0,b)$ and
is continuous at $b$, it is Lipschitz on the segment
$[0,b]$. However, the variation of $h$ between $a_n$ and $b_n$ is
\begin{eqnarray*}
\bigvee\limits_{a_n}^{b_n}(h)
&=&
\sum\limits_{i=1}^{2k_n} |h(a_n+id_n) - h(a_n+(i-1)d_n)|\\
&=&
\sum_{i=1}^{2k_n} |f(x_n)-f(y_n)| \geq n \sum_{i=1}^{2k_n}d_n =2n\,d_nk_n\sim \frac{2}{n},
\end{eqnarray*}
and, therefore,
$\bigvee\limits_{0}^{b}(h)\geq \sum_{n=1}^{\infty}\,\bigvee \limits_{a_n}^{b_n}(h)=\infty$.
Thus $h$ is not absolutely continuous on $I = [0, b]$.
\qed

\medskip
\noindent
\textbf{Remark 1}. Theorem \ref{solution} is false for
convex sets $P$ in linear metric spaces. Let, for example, $0<p<1$
and $P=[0,1]\subset (\mathbb R,|\cdot|_p)$, with the distance
$|x-y|_p:=|x-y|^p$. It is not difficult to verify that every
absolutely continuous function $g:I\to P$ is constant. Therefore,
the set $P$ and an arbitrary mapping $f$ satisfy the conditions
$(ii)$ and $(iii)$. But, for example, the function $f:P\to\mathbb
R$, defined by $f(x)=x^{p/2}$, is not Lipschitz because
$$
\frac{|x^{p/2}-0|}{|x-0|_p}=x^{-p/2}\to\infty\;\;\; \mbox{as}\;\;\;x\to 0^+.$$

\begin{cor}\label{2variables}
Let $f:Q\to\mathbb R$ be a function. The following conditions are
equivalent:
\begin{itemize}
\item[$(i)$] $f$ is Lipschitz on $Q$.

\item[$(ii)$] For every absolutely continuous functions
$g_1,g_2:I\to I$ the superposition $f(g_1(x),g_2(x))$ is
absolutely continuous.

\item[$(iii)$] For every Lipschitz functions $g_1,g_2:I\to I$ the
superposition $f(g_1(x),g_2(x))$ is absolutely continuous.
\end{itemize}
\end{cor}

{\it Proof.} Since the mapping $g=(g_1,g_2):I\to Q$ is Lipschitz
(absolutely continuous) if and only if the functions $g_1, g_2$
are Lipschitz (absolutely continuous), Theorem~\ref{solution}
implies Corollary \ref{2variables}. \qed \bigskip

\section{The diagonal case}
By putting  $g_1(x)=g_2(x)=x$ in Eidelheit's problem, we obtain
the question on absolute continuity for the diagonal of a
separately absolutely continuous function. In this section we give
a negative answer to a stronger version of this question.

\begin{theo}\label{diagonalTheo}
Let $(u_n)^{\infty}_{n=1}$ be a sequence of reals $u_n>0$ such
that $\sum_{n=1}^\infty u_n<\infty$. Then there exists a separately
Lipschitz function $f:Q \to\mathbb R$ such that the Lebesgue measure
\begin{equation}\label{Mykh1}
\lambda(\{z\in Q:f'_x(z)\ne 0\,\,\,{\rm or}\,\,\, f'_y(z)\ne 0\})\le \sum\nolimits_{n=1}^{\infty} u_n,
\end{equation}
\begin{equation}\label{Mykh2}
\lambda(\{z\in Q:|f'_x(z)|\geq 2^n\})\le \sum\nolimits_{k=n}^{\infty} u_k
\end{equation}
and
\begin{equation}\label{Mykh3}
\lambda(\{z\in Q:|f'_y(z)|\ge 2^n\})\le \sum\nolimits_{k=n}^{\infty} u_k
\end{equation}
for every $n\in\mathbb N$, and the function $h(x)=f(x,x)$ has unbounded variation on $I$.
\end{theo}

{\it Proof.} For every $n\in\mathbb N$ let
$$
I_n=[\frac{1}{2^n},\frac{1}{2^{n-1}}], ~ k_n=[\frac{1}{4^nu_n}]+1~~ {\rm and} ~~ d_n=\frac{1}{2^nk_n}.
$$
Note that $k_n\, d_n^2\leq u_n$. Moreover, for every $n\in\mathbb
N$ and $1\leq i\leq k_n$ let
$$
I_{n,i}=[a_{n,i},
b_{n,i}]=\left[\frac{1}{2^n}+(i-1)d_n,\frac{1}{2^{n}}+id_n\right],
~ c_{n,i}=\frac{a_{n,i}+b_{n,i}}{2}, ~ Q_{n,i}=I_{n,i}\times
I_{n,i}
$$
and choose a continuous separately Lipschitz function
$\varphi_{n,i}:Q_{n,i} \to\mathbb R$ so that
\begin{itemize}
\item[$(1)$] $\varphi_{n,i}(x,y)=0$ if $(x,y)\not\in\mathrm{int}\,Q_{n,i}$,

\item[$(2)$] $\varphi_{n,i}(c_{n,i},c_{n,i})=\frac{1}{2 k_n}$,

\item[$(3)$]  $\varphi_{n,i}$ is linear on every segment connecting a
boundary point of $Q_{n,i}$ with $(c_{n,i},c_{n,i})$.
\end{itemize}
Next, put
$$
f(x,y)=\sum_{n=1}^{\infty}\sum_{i=1}^{k_n}\varphi_{n,i}(x,y).
$$
Since all segments $(a_{n,i},b_{n,i})$ are disjoint, for every $x_0,y_0\in(a_{n,i},b_{n,i})$
and $x,y\in I$ we have that $f(x_0,y)=\varphi_{n,i}(x_0,y)$ and $f(x,y_0)=\varphi_{n,i}(x,y_0)$.
Therefore, $f$ is separately Lipschitz.

\medskip
For every $n\in\mathbb N$ consider
$$
A_n=\{z\in Q:|f'_x(z)|\geq 2^n\} ~ {\rm and} ~~ B_n=\{z\in Q:|f'_y(z)|\geq 2^n\}.
$$
Note that
$|(\varphi_{n,i})'_x(z)|\leq \frac{1}{d_nk_n}=2^n$ almost everywhere
on $Q_{n,i}$. Therefore
$$
A_n\subseteq \bigcup\nolimits_{m=n}^\infty \bigcup\nolimits_{i=1}^{k_m} Q_{m,i},
$$
and thus
$$
\lambda(A_n)\leq \sum_{m=n}^\infty \sum_{i=1}^{k_m} \lambda(Q_{m,i})=
\sum_{m=n}^\infty k_md_m^2\leq\sum_{m=n}^{\infty} u_m.
$$
Similarly we find that $\lambda(B_n)\leq\sum_{m=n}^{\infty} u_m$. Since
$$
C:= \{z\in Q:f'_x(z)\ne 0\,\,\,{\rm or}\,\,\, f'_y(z)\ne 0\}\subseteq
\bigcup\nolimits_{n=1}^\infty \bigcup\nolimits_{i=1}^{k_n} Q_{n,i}
$$
it follows that
$$
\lambda(C)\leq\sum_{n=1}^\infty \sum_{i=1}^{k_n} \lambda(Q_{n,i})=\sum_{n=1}^\infty u_n.
$$
We only need to show that the function $h(x)=f(x,x)$ has unbounded variation on $I$. We have that
\begin{eqnarray*}
\bigvee\limits_{0}^{1}(h)
&\geq&
\sum_{n=1}^{\infty} \sum_{i=1}^{k_n} \left(| h(c_{n,i} - h(a_{n, i}))| + | h(b_{n,i} - h(c_{n, i}))| \right) \\
&=& \sum_{n=1}^{\infty} \sum_{i=1}^{k_n} \left(\frac{1}{2k_n} +
\frac{1}{2 k_n}\right) = \infty.
\end{eqnarray*}
The proof is complete. \qed

\medskip
A Banach space $E$ of classes of measurable functions $f:Q\to\mathbb R$ is called
{\it rearrangement invariant $($r.i.$)$ Banach function space} or {\it symmetric space}
(over $Q$ with the Lebesgue measure $\lambda$ such that $\lambda(Q) = 1$), if it satisfies
the following conditions (see \cite[p. 59]{BS}, \cite[p. 90]{KPS} and \cite[pp. 114-119]{LT2}):
\begin{itemize}
\item[$(1)$] if $|f(z)|\leq|g(z)|$ for $\lambda$-almost every $z\in Q$, $f$
measurable on $Q$ and $g\in E$, then $g\in E$ and $\|f\|_E\leq \|g\|_E$;

\item[$(2)$] if $f$ and $g$ are equimeasurable, that is,
$\lambda(\{z\in Q: |f(z)| > \alpha\}) = \lambda(\{z\in Q: |g(z)| >
\alpha\})$ for every $\alpha \geq 0$ and $g\in E$, then $f\in E$
and $\|f\|_E= \|g\|_E$.

\end{itemize}

\medskip
\noindent
Note that for any r.i. space $E$ on $Q$ we have continuous embeddings
$$
L_{\infty}(Q) \subset E \subset L_1(Q) \,\, {\rm ~with} \,\, \| f\|_{L_1(Q)} \leq \frac{\| f\|_E}{\|\chi_Q\|_E} \leq \| f\|_{L_{\infty}(Q)} \, \, {\rm ~for ~all } \,\, f \in L_{\infty}(Q).
$$
Moreover, since $\lambda(Q) = 1$ we can have as the definition of equimeasurability
of $f$ and $g$ in (2) the equality $\lambda(\{z\in Q: |f(z)| \geq \alpha\}) = \lambda(\{z\in Q: |g(z)| \geq
\alpha\})$ for every $\alpha > 0$.

\bigskip
The following lemma is well known (see e.g. \cite[p. 2]{Brav}, \cite[p. 98]{KPS}):

\begin{lem}\label{r.i.}
Let $E$ be a r.i. Banach function space on $Q$, $g\in E$ and $f: Q\to\mathbb R$
be a measurable function such that
$$
\lambda(\{z\in Q: |f(z)|\geq \alpha\}) \leq \lambda(\{z\in Q: |g(z)|\geq \alpha\})
$$
for every $\alpha \geq 0$. Then $f\in E$ and $\|f\|_E\leq \|g\|_E$.
\end{lem}

\medskip
\begin{lem}\label{seq}
Let $(v_n)^{\infty}_{n=1}$ be a decreasing sequence of reals $v_n>0$.
Then there exists a sequence
$(u_n)^{\infty}_{n=1}$ of reals $u_n>0$ such that $\sum_{k\geq n}u_k\le v_{n}$
for every $n$.
\end{lem}

{\it Proof.} It is sufficient to take
$u_n=\frac{v_n}{n}-\frac{v_{n+1}}{n+1}$ for every $n$. \qed

\begin{cor}\label{diagonalCor1}
Let $\{E_s\}_{s\in S}$ be a family of r.i. Banach function spaces $E_s$ on $Q$ such that
$(\bigcap\nolimits_{s\in S}E_s)\setminus L_\infty(Q)\ne \O$. Then there
exists a separately Lipschitz function $f:Q \to\mathbb R$
such that $f'_x, f'_y\in \bigcap\nolimits_{s\in S}E_s$ and the function
$h(x)=f(x,x)$ has unbounded variation.
\end{cor}

{\it Proof.} First, let us take a function $g\in (\bigcap\nolimits_{s\in S}E_s)\setminus L_\infty(Q)$
and for every $n\in\mathbb N$ put
$$
v_n=\lambda(\{z\in Q:|g(z)|\geq 2^{n+1}\}).
$$
By using Lemma \ref{seq} we can choose a sequence $(u_n)^{\infty}_{n=1}$ of reals
$u_n>0$ such that $\sum_{k=n}^{\infty} u_k\leq v_{n}$ for every $n\in\mathbb N$.

Theorem \ref{diagonalTheo} implies that there exists a separately
Lipschitz function $f:Q \to\mathbb R$ such that the conditions
(\ref{Mykh1})--(\ref{Mykh3}) from Theorem \ref{diagonalTheo} are
satisfied for every $n\in\mathbb N$, and the function
$h(x)=f(x,x)$ has unbounded variation.

We only need to show that $f'_x, f'_y\in \bigcap\nolimits_{s\in S}E_s$. For fixed $\alpha>0$
let
$$
A_\alpha=\{z\in Q: |f'_x(z)|\ge \alpha\} ~~{\rm and} ~ B_\alpha=\{z\in Q: |g(z)|\ge \alpha\}.
$$
If $\alpha\in(0,4]$, then
\begin{eqnarray*}
\lambda(A_\alpha)
&\leq&
\lambda(\{z\in Q:f'_x(z)\ne 0\}) \leq \sum_{n=1}^{\infty} u_n\le v_1\\
&=&
\lambda(\{z\in Q:|g(z)|\geq 4\})\leq \lambda(B_\alpha).
\end{eqnarray*}
If $\alpha>4$, then choosing $n\in\mathbb N$ with
$2^n<\alpha\le 2^{n+1}$ we have that
\begin{eqnarray*}
\lambda(A_\alpha)
&\leq&
\lambda(\{z\in Q:|f'_x(z)|\ge 2^n\})\le \sum_{k=n}^{\infty} u_k \leq v_n\\
&=&
\lambda(\{z\in Q:|g(z)|\ge 2^{n+1}\})\le \lambda(B_\alpha).
\end{eqnarray*}
Thus $\lambda(A_\alpha)\le \lambda(B_\alpha)$ and
$f'_x\in \bigcap_{s\in S}E_s$ by Lemma \ref{r.i.}.
Similarly we can prove that $f'_y\in \bigcap_{s\in S}E_s$.
The proof is complete. \qed

\begin{cor}\label{diagonalCor2}
There exists a separately Lipschitz function $f:Q \to\mathbb R$
such that $\mathop{\int\!\!\!\int}\limits_{Q}|f'_x|^p
\,dxdy<\infty$ and $\mathop{\int\!\!\!\int}\limits_{Q}|f'_y|^p
\,dxdy<\infty$ for every $p>1$, and the function $h(x)=f(x,x)$ has
unbounded variation.
\end{cor}

{\it Proof.} Note that for every $p>1$ the space $L_p(Q)$ is r.i. space and the
function $g(x,y)=\ln x$ belongs to $\bigcap_{p>1}L_p(Q)\setminus L_\infty(Q)$.
\qed \bigskip

\noindent{\bf Remark 2.} Note that $g(x_2)-g(x_1)=\mathop{\int}\limits_{x_1}^{x_2}g\prime
\,dx$ for every absolutely continuous function $g:I\to\mathbb R$.
If  the partial derivatives $f'_x$ and $f'_y$ of a separately
absolutely continuous function $f:Q\to\mathbb R$ are such that $|f'_x|\leq C$
and $|f'_y|\leq C$ almost everywhere on $Q$,
then Fubini theorem implies that $f$ is separately Lipschitz with
the constant $C$ and, therefore, $f$ is jointly Lipschitz with the constant $C$
with respect to the sum-distance on $Q$, in particular, the restriction of $f$
on any straight line is Lipschitz.

\section{``Embeddings of Banach spaces'' approach}

We show how, using Theorems A and 2, one can give an answer (in a
classical Banach style) to the Eidelheit question. Moreover, we
obtain, as a byproduct, stronger results, which are not evident
under the function theory approach. \smallskip

Our approach is based on the following well known notion: A
bounded linear operator $T$ from a topological vector space $X$
into a topological vector space $Y$ is called {\it strictly
singular} if there exists no infinite dimensional subspace
$Z\subset X$ such that $T|_Z$ is an isomorphism. The operator $T$
from a Banach space $X$ into a Banach space $Y$ is called {\it
superstrictly singular} ({\it SSS} for short) if there does not
exist a number $c>0$ and a sequence of subspaces $E_{n} \subset
X$, $\dim E_{n}=n$, such that $\|Tx\| \ge c \|x\|$ for all $x$ in
$\bigcup_n E_n.$ Obviously, each compact operator is SSS, each SSS
operator is strictly singular and $T$ is SSS if it is SSS on a
finite codimensional closed subspace (cf. \cite{Plichko}).

A Rudin's version \cite[Th. 5.2]{Rudin0} of the well-known
Grothendieck's result says that the natural (noncompact) embedding
$\mathcal{I}: L_{\infty}\hookrightarrow L_p$, $p\ge 1$, is SSS. On
the other hand, generalizing  the Grothendieck's result, Novikov
\cite[Th. 1]{Novikov} has proved that the natural embedding
$\mathcal{I}: L_{\infty}\hookrightarrow E$ is strictly singular for
every rearrangement invariant (r.i.) space $E\ne
L_{\infty}$. The following result contains both the Rudin and
Novikov theorems.

\begin{theo}\label{SSS}
Let $E$ be a r.i. Banach function space on $I = [0, 1]$ different from $L_{\infty}(I)$.
Then the natural embedding $\mathcal{I}:L_{\infty}\hookrightarrow
E$ is SSS.
\end{theo}

In the proof we will use the following well-known lemma (see e.g. \cite[Lemma 3.3]{Popov}).

\medskip
\begin{lem}\label{intersection}
Let $b>0$. Then for every $k\in \mathbb N, k\geq 2$ there exists $n = n(b, k) \in \mathbb N$ such that
for any collection of measurable subsets $A_i\subset I,\; i=1,\dots,n$ with the
Lebesgue measure $\lambda (A_i)>b$ there is a subcollection $(A_{i_j})_{j=1}^k$ with
\begin{equation}\label{Popov}
\lambda\left(\bigcap\nolimits_{j=1}^k A_{i_j}\right)>0\;.
\end{equation}
\end{lem} \medskip

{\it Proof of Theorem} \ref{SSS}. Let us suppose on the contrary, that there exist
$\varepsilon >0$ and  $n$-dimensional subspaces $E_n\subset L_{\infty}$
such that for every $n$ and $f\in E_n$
\begin{equation}\label{Novikov}
\varepsilon \|f\|_{L_{\infty}}\le \|f\|_E\;.
\end{equation}
Since $E\ne L_{\infty}$, there are $a$ and $b$ such that for every
$f\in L_{\infty}$ with $\|f\|_{L_{\infty}}=1$ and $\|f\|_E\ge
\varepsilon,$
$$
\lambda(\{x:|f(x)|>a\}) >b.
$$
Then, by Lemma~\ref{intersection}, there exists $c>0$ such that for every
$k$ with the property (\ref{Popov}) there is $n$ so that for any
elements $(f_i)_1^k$ with the property (\ref{Novikov}) we have that
$$
\frac{1}{k}\left\|\sum\nolimits_{i=1}^kf_i\right\|_E\ge c\;.
$$

\newpage
\noindent
Take an orthogonal (with respect to natural inner product) basis
$(f_i^n)_{i=1}^n$ of $E_n$ with $\|f_i^n\|_{L_{\infty}}=1$. Then
$$
\left\|\sum\nolimits_{i=1}^nf_i^n\right\|_{L_2}\le n^{1/2}\;.
$$
Put
$$
\sigma(n,\delta )=\left\{x\in (0,1):\frac{1}{n}\left|\sum\nolimits^n_{i=1}f_i^n(t)\right|>\delta\right \}\;.
$$
Hence, for every $\delta>0$ the measure $\;\lambda(\sigma(n,\delta))\to 0$ as $n\to \infty$ not
depending on the form of $(f_i)$ (see e.g. \cite[ p. 160]{LT2}). But
$$
\frac{1}{n}\left\|\sum\nolimits_{i=1}^nf_i^n \right\|_E\le \delta
+\|\chi_{\sigma(n,\delta)}\|_E\;.
$$
Since $E$ is different from $L_{\infty}\,$ we have $\|\chi_{A_n}\|_E\to 0$ provided
$\lambda (A_n)\to 0$ (see e.g. \cite[p. 118]{LT2}). Thus,
$$
\frac{1}{n}\left\|\sum\nolimits_{i=1}^nf_i^n\right\|_E\to
0\;\;\mathrm{as}\;\;n\to \infty,
$$
and we have a contradiction.
 \qed

\medskip
\noindent
{\bf Remark 3.} Of course, Theorem \ref{SSS} is valid for r.i. Banach function spaces
on any subset $A\subset \mathbb{R}^n$ of positive finite Lebesgue measure.

\medskip
Let us look at the integral condition in the one variable version of Eidelheit's problem.
It seems that he means the existence of the derivative almost everywhere on
$I$, i.e. Eidelheit had generalized derivatives. The corresponding one variable
space was considered as far back as by Banach. Namely, in \cite[pp. 134,
167]{Banach} he introduced, in particular, a space of absolutely
continuous functions on $I$ with derivative in $L^p(I)$. On this space one can
introduce the norm $\|f\|=|m_f|+\|f'\|_{L_p}$ (this is just one of equivalent forms).
Banach noted that this space is in fact complete.

Given an arbitrary r.i. Banach function space $E$ on $I$, one can define
the Beppo Levi space $BL^1_E(I)$ (why this space is named after Beppo
Levi, we will explain below) of absolutely continuous functions
$f$ for which $f'\in E$ with the natural norm (this is just one of the equivalent norms):
$$\|f\|=|m_f|+\|f'\|_E\,$$ where $m_f$ is the mean of the function
$f$ on $I$.

\noindent
{\bf Remark 4.} Let $Y$ be the one codimensional subspace
of $BL^1_{E}(I)$ consisting of all $f$ with $m_f=0$ and let
$D':Y\to E$ be the generalized derivative operator. Obviously,
this operator is an onto isometry, so $BL^1_{E}(I)$ is complete.

\begin{cor}\label{SSSCor4}
Let $E$ be a r.i. Banach function space on $I$ different from $L_{\infty}$.
Then the natural embedding $\mathcal{J}:\Lip_1(I)\hookrightarrow
BL^1_E(I)$ is non-compact but SSS.
\end{cor}

{\it Proof.} Take the Rademacher functions $(r_n)$ and put
$f_n(x)=\int_0^x r_n(t)\,dt\,,\;x\in I$ (the Schauder functions).
Then $f_n\in\Lip_1(I)$ and $\|f_n\|=1$ for every $n$. On the other
hand, $\|\mathcal{J}f_n-\mathcal{J}f_m\|\ge 1$ for any $n\ne m$.
Thus, mapping $\mathcal{J}$ is not compact.

Denote by $X$ the subspace of $\Lip_1(I)$ consisting of all $f$
with $m_f=0$. Then the generalized derivation operator $D:X\to
L_{\infty}$ defined by $Df=f'$ is bounded. Let $Y$ and $D'$ be
from Remark 4. Then from the following operator diagram

\begin{figure}[!h]
\setlength{\unitlength}{1.0cm}
\small
\centering
\begin{picture}(6,3)(0,0)
\put(2.2,2.5){\vector(1,0){1.7}}
\put(4.15,2.5){\makebox(0,0){$E$}}
\put(4.1,0.7){\vector(0,1){1.5}}
\put(4.15,0.4){\makebox(0,0){$Y$}}
\put(1.6,0.7){\vector(0,1){1.5}}
\put(2.0,0.5){\vector(1,0){1.8}}

\put(1.55,2.5){\makebox(0,0){$L_\infty(I)$}}
\put(1.55,0.4){\makebox(0,0){$X$}}

\put(3.0,2.75){\makebox(0,0){$\mathcal{I}$}}
\put(3.0,0.25){\makebox(0,0){$\mathcal{J}$}}

\put(1.2,1.5){\makebox(0,0){$D$}}
\put(4.6,1.5){\makebox(0,0){$D^{\prime}$}}
\end{picture}
\label{fig:opdiagram}
\end{figure}

\noindent
we obtain that $\mathcal{J}=(D')^{-1}\mathcal{I}D$.  Since, by
Theorem \ref{SSS}, $\mathcal{I}$ is SSS, it follows that the restriction
$\mathcal{J}|_X$ is SSS, and so $\mathcal{J}$ is SSS. The proof is complete. \qed
\bigskip

Let us consider a sequence of r.i. Banach function spaces $(E_k)$ different from
$L_{\infty}(I)$ and the topological vector space $F
=\bigcap_{k}BL^1_{E_k}(I)$ with the fundamental neighborhood of
$0$ formed by the unit balls of the spaces $BL^1_{E_k}(I)$,
$k=1,2,\dots$ It is easy to verify that $F$ is a Fr\'echet space.

\begin{cor}\label{SSCor5}
The natural embedding $\mathcal{J}:\Lip_1(I)\hookrightarrow F$ is
strictly singular.
\end{cor}

{\it Proof.} Let $Z\subset \Lip_1(I)$ be an infinite dimensional
subspace and $B_Z$ be its unit ball. If $\mathcal{J}|_Z$ is an
isomorphism, then $\mathcal{J}(B_Z)$ is a neighborhood of $0$ in
$\mathcal{J}(Z)\subset F$. Hence there exists $k$ such that
$\mathcal{J}(B_Z)$ is a neighborhood of $0$ in
$\mathcal{J}(Z)\subset BL^1_{E_k}(I)$. This contradicts to
Corollary \ref{SSSCor4} and the proof follows. \qed
\bigskip

The next corollary generalizes the solution of the one variable
version of Eidelheit's problem.

\begin{cor}\label{Cor6}
Every infinite dimensional closed subspace of $F$ contains an
absolutely continuous function $f$ for which there is a Lipschitz
function $g$ with the non-absolutely continuous superposition
$f\circ g$.
\end{cor}

{\it Proof.} Indeed, by Corollary \ref{SSCor5}, every infinite
dimensional closed subspace of $F$ contains an (absolutely
continuous) function $f$ which does not belong to $\Lip_1(I)$. By
Theorem A, there is a Lipschitz function $g$ with the
non-absolutely continuous superposition $f\circ g$. \qed \bigskip

Before presenting the abstract versions of the two variable results and
the diagonal Eidelheit problems let us start with a short historical excursion.
The conditions in Eidelheit problem mean that $f\in BL_p^1(Q)$, where
the (Beppo Levi) space $BL_p^1(Q)$ consists of all functions
$f:Q \rightarrow \mathbb R$, which are absolutely continuous in each
variable and whose (classical generalized) first order partial
derivatives lie in $L_p(Q)$. For $p=2$ and three variables a
similar space was considered by O. Nikodym \cite{Nikodym}. This
space is called {\it Beppo Levi space} since functions in this
class were studied as far back as 1906 by Beppo Levi (for $p=2$),
and later by Tonelli (for $p=1$ and $p\ge 2$) in the minimization
of variational integrals. The name Beppo Levi space was introduced by
Nikodym \cite{Nikodym} for $p = 2$ in 1933 and in general by
Deny-Lions \cite{DenyLions} in 1953. On the other hand, in 1936
Sobolev \cite{Sobolev} developed the so called {\it Sobolev space}
$W^1_p(Q)$ as a space of all $f\in L_p(Q)$
whose generalized (distributional) derivatives belong to $L_p(Q)$.
Surprisingly, we have that $BL_p^1(\Omega) = W^1_p(\Omega)$ even for
subsets $\Omega \subset \mathbb R^n$ (cf. \cite[Th. 1, p. 8]{Mazya}
and \cite[Th. 2.1.4]{Ziemer}) but we must explain in which sense
since $BL_p(Q)$ is a space of functions and $W^1_p(Q)$ a space of
equivalence classes of functions which differ on sets of measure
zero. More correctly, we mean that any function in $BL_p(Q)$ lies
in (an equivalence class in) $W^1_p(Q)$, while every element of
$W^1_p(Q)$ has a representative in $BL_p(Q)$. \vspace{3mm}

Denote by $\Lip_1(Q)$ the space of Lipschitz functions on $Q$ with
the norm
$$
\|f\|=|m_f|+\sup_{z,z'\in Q,\;z\ne z'}\frac{|f(z)-f(z')|}{d(z,z')}\;,
$$
where $m_f=\mathop{\int\!\!\!\int}_Q f(x,y)\,dxdy$ and $d$ denotes
the Euclidean distance in $\mathbb{R}^2$.

\medskip
Let $E$ be a r.i. Banach function space on $Q$. Denote by $BP^1_E(Q)$
the space of functions $f$ on $Q$, which are absolutely continuous with respect
to each variable for almost all other variables, and whose
generalized partial derivatives $f'_x, f'_y\in E$ with the natural norm
$$
\|f\|=|m_f|+\|f'_x\|_E+\|f'_y\|_E\,.
$$
Similar spaces were considered by Deny and Lions \cite{DenyLions}.
They mean the derivatives in the sense of generalized functions.
Deny and Lions have proved that these spaces are complete.

\medskip

\noindent{\bf Remark 5.} Let $Y$ be the one codimensional subspace
of $E$ consisting of all $f$ with $m_f=0$. Let $D': Y\to E\times
E$ be defined by $D'(f)=(f'_x,f'_y)$. This operator is an into
isomorphism. So, $BL^1_{E}(Q)$ is isomorphic to a subspace of $E$.
Is it isomorphic to $E$? For the spaces $BL^1_{p}(Q)$ the answer
is well-known \cite{PelczynskiWojciechowski}.

\begin{cor}\label{SSSofBL}
Let $E$ be a Banach r.i. space on $Q$ different from
$L_{\infty}(Q)$. Then the natural embedding
$\mathcal{J}:\Lip_1(Q)\hookrightarrow BL^1_{E}(Q)$ is SSS.
\end{cor}

{\it Proof.} Denote by $X$ the subspace of $\Lip_1(Q)$ consisting
of all $f$ with $m_f=0$. Then the generalized derivation operator
$D:X\to L_{\infty}(Q)\times L_{\infty}(Q)$, $Df=(f'_x,f'_y)$ is
bounded. Let $\mathcal{I}: L_{\infty}(Q)\times
L_{\infty}(Q)\hookrightarrow E\times E$ be the natural embedding.
By Remark 1, it is SSS.

\noindent
Let $Y$ and $D'$ be from Remark 4. Then from the following diagram

\begin{figure}[!h]
\setlength{\unitlength}{1.0cm}
\small
\centering
\begin{picture}(6,3)(0,0)
\put(3.2,2.5){\vector(1,0){1.8}}
\put(5.75,2.5){\makebox(0,0){$E \times E$}}
\put(5.75,0.7){\vector(0,1){1.5}}
\put(5.75,0.4){\makebox(0,0){$Y$}}
\put(1.6,0.7){\vector(0,1){1.5}}
\put(2.0,0.5){\vector(1,0){3.4}}

\put(1.55,2.5){\makebox(0,0){$L_\infty(Q) \times L_\infty(Q)$}}
\put(1.55,0.4){\makebox(0,0){$X$}}

\put(4.0,2.75){\makebox(0,0){$\mathcal{I}$}}
\put(4.0,0.25){\makebox(0,0){$\mathcal{J}$}}

\put(1.2,1.5){\makebox(0,0){$D$}}
\put(6.25,1.5){\makebox(0,0){$D^{\prime}$}}
\end{picture}
\label{fig:opdiagram}
\end{figure}
\noindent
we have that $\mathcal{J}=(D')^{-1}\mathcal{I}D$. Hence,
$\mathcal{J}|_X$ is SSS, and so $\mathcal{J}$ is SSS.
\qed

\medskip
Let us consider a sequence of r.i. Banach function spaces $(E_k)$
on $Q$, different from $L_{\infty}(Q)$, and the topological vector
space $F =\bigcap_{k}BL^1_{E_k}(Q)$ with the fundamental
neighborhood of $0$ formed by unit balls of the spaces
$BL^1_{E_k}(Q)$, $k=1,2,\dots$ It is easy to verify that $F$ is a
Fr\'echet space. The proof of the following corollary is the same
as in Corollary \ref{SSCor5}.

\begin{cor}\label{Cor8}
The natural embedding $\mathcal{J}:\Lip_1(Q)\hookrightarrow F$ is
strictly singular.
\end{cor}

The next corollary generalizes the solution of Eidelheit's problem.

\begin{cor}\label{Cor9}
Every infinite dimensional closed subspace of $F$ contains an
$($absolutely continuous$)$ function $f$ for which there are
Lipschitz functions $g_1,g_2$ with the non-absolutely continuous
superposition $f(g_1(t),g_2(t))$.
\end{cor}

{\it Proof.} Indeed, by Corollary \ref{Cor8}, every infinite
dimensional closed subspace of $F$ contains an (absolutely
continuous) function $f$ which does not belong to $\Lip_1(Q)$. By
Corollary 1, there are Lipschitz functions $g_1,g_2$ with the
non-absolutely continuous superposition $f(g_1(t),g_2(t))$. \qed
\bigskip

Let us now consider the functional analytic meaning of
Theorem~\ref{diagonalTheo}. Denote by $Z$ the ``diagonal''
subspace of $BL^1_{E}(Q)$ consisting of functions $f(x,y)\in
BL^1_{E}(Q)$ for which $f(x+\lambda x,x-\lambda x)=f(x,x)$ ,
$\,x\in I$ , $\,\lambda \in \mathbb{R}$. From Corollary \ref{SSSofBL}
we have then immediately that:

\begin{cor}\label{Cor10}
There exists $f\in Z$ such that $f\not\in\Lip_1(Q)$.
\end{cor}\medskip

Note that Corollary \ref{diagonalCor2} is stronger than Corollary
\ref{Cor10} since in Corollary \ref{diagonalCor2} $f$ is a separately
Lipschitz function.


\vspace{3mm}

\noindent
{\footnotesize Department of Mathematics, Lule\r{a} University of Technology\\
SE-971 87 Lule\r{a}, Sweden} ~{\it e-mail:} {\tt lech@sm.luth.se} \\

\vspace{-3mm}

\noindent
{\footnotesize Department of Applied Mathematics, Chernivtsi National University\\
Kotsyubyn'skoho 2, 58012 Chernivtsi, Ukraine} ~{\it e-mail:} {\tt mathan@chnu.cv.ua} \\

\vspace{-3mm}

\noindent
{\footnotesize Institute of Mathematics, Cracow University of Technology,\\
Warszawska 24, 31-155 Krak\'ow, Poland}  {\it e-mail:} {\tt aplichko@usk.pk.edu.pl} \\

\end{document}